\documentclass[10pt]{amsart}

\usepackage{amssymb}
\usepackage{amsmath}
\usepackage{amsthm}

\input xy
\xyoption{all}

\title[Universal graph operator algebras: constructions and uniqueness]
{Explicit construction and uniqueness for universal operator algebras of directed graphs}

\author{Benton L. Duncan}

\address{Department of Mathematics\\
300 Minard Hall\\
North Dakota State University\\
Fargo, ND  58105-5075}

\email{Benton.Duncan@ndsu.nodak.edu}



\thanks{Part of this work was supported by a Department of Education Grant \#P200A030193}

\begin{document}

\theoremstyle{plain}
\newtheorem{thm}{Theorem}
\newtheorem{lem}{Lemma}
\newtheorem{prop}{Proposition}
\newtheorem{cor}{Corollary}

\theoremstyle{definition}
\newtheorem{dfn}{Definition}
\newtheorem*{construction}{Construction}[section]
\newtheorem*{example}{Example}[section]

\theoremstyle{remark}
\newtheorem*{question}{Question}
\newtheorem*{acknowledgement}{Acknowledgements}
\newtheorem{rmk}{Remark}

\begin{abstract} Given a directed graph, there exists a universal
operator algebra and universal $C^*$-algebra associated to the
directed graph.  In this paper we give intrinsic constructions of
these objects.  We provide an explicit construction for the
maximal $C^*$-algebra of an operator algebra. We also discuss
uniqueness of the universal algebras for finite graphs, showing
that for finite graphs the graph is an isomorphism invariant for
the universal operator algebra of a directed graph.  We show that
the underlying undirected graph is a Banach algebra isomorphism
invariant for the universal $C^*$-algebra of a directed graph.
\end{abstract}

\maketitle

There has been significant work in the study of operator algebras
associated to combinatorial objects (e.g. groups, semigroups, and
graphs).  We have continued this study in \cite{Duncan:2004} where
the universal operator algebra of a directed graph and the
universal $C^*$-algebra of a directed graph were introduced and
described. The aim of this paper is twofold: first we refine the
construction of the universal operator algebras of directed
graphs, then we discuss invariants of the universal algebras of
finite directed graphs.

First we use ideas from \cite{Bl-Paul:1991} to define intrinsic
norms on $OA(Q)$ the universal operator algebra of a directed
graph. This allows a more concrete construction than was given in
\cite{Duncan:2004}. We also describe a construction of the maximal
$C^*$-envelope of an operator algebra, see \cite{Bl:1999}. This
construction is defined instrinisically using the free product
operator algebra construction of Blecher and Paulsen
\cite{Bl-Paul:1991}.  This suggests that the maximal
$C^*$-envelope is not as mysterious as is presumed.  In fact
having a canonical construction should allow a more detailed study
of the maximal operator algebra of a directed graph in particular
cases.

Kribs and Power show in \cite{Kribs-Power:2003a} that the graph is
a complete unitary invariant for the Toeplitz quiver algebra of a
directed graph. Recent work on these Toeplitz quiver algebras by
Katsoulis and Kribs, \cite{Kat-Kribs:2003} and by Solel
\cite{Sol:2003}, has demonstrated that the graph is a complete
isomorphism invariant for these algebras.  In this paper we extend
the techniques of \cite{Kat-Kribs:2003} to show that for finite
graphs the graph is a complete isomorphism invariant for $OA(Q)$.
This fact is perhaps not surprising, although the technique
requires more subtlety than in \cite{Kat-Kribs:2003} and
\cite{Sol:2003}.  For the universal $C^*$-algebra of a finite
directed graph we are able to show that the underlying directed
graph is an isomorphism invariant for the algebra. This is very
surprising since the Cuntz-Krieger algebra of a directed graph is
not classified by the graph.

Before proceeding, we would like to emphasize a difference between
the operator algebras in the present paper and those defined in
\cite{Muhly:1997}. When we construct the universal operator
algebra of a directed graph we consider representations which send
vertices to projections.  \emph{We do not assume that the
projections are orthogonal}, as was implicit in \cite{Muhly:1997}.
This difference provides examples which differ significantly from
the Toeplitz quiver algebras.

We remind the reader of some definitions which can be found in
\cite{Duncan:2004}.  If $Q$ is a directed graph we will let $V$
and $E$ be the vertex and edge sets, respectively.  We let $W(Q)$
be the set of finite words in $E \cup V$ and make $W(Q)$ a
semigroup by defining multiplication via concatenation.

We would also like to clarify the construction of the universal
operator algebra of a directed graph from \cite{Duncan:2004}.  It
was implicit in the construction that the set $ w(Q)$ was the set
of all finite words in $V(Q) \cup E(Q)$ subject to the relations $
r(e) e = e = e s(e)$.  Here $r: E\cup V \rightarrow V$ and $ s: E
\cup V \rightarrow V$ are the range and source map extended to $E
\cup V$ by defining $r(v) = v = s(v)$ for all $ v \in V$.

We point out the connection between Section 2 and the results in
\cite{Palm:2001}.  There the Gelfand-Naimark seminorm is defined
on a Banach $*$-algebra.  It turns out that the seminorm we define
in Section 2 is equal to the Gelfand-Naimark seminorm. This yields
a more direct approach the material in Section 2.  We have used
both methods, we concretely define the seminorm to emphasize that
the norm is intrinsic, and then we reference Palmer's work for
completeness.

We now emphasize an important fact from the latter parts of this
paper.  \emph{In sections 4, 5 and 6 we restrict our attention to
finite graphs.}

\section{Intrinsic norms for $OA(Q)$}

We begin by adapting a technique of \cite{Bl-Paul:1991} to provide
intrinsic norms on $OA(Q)$.  Recall that for a monoid $M$ an
intrinsic norm is described in \cite{Bl-Paul:1991} for the
universal operator algebra associated to a monoid, denoted $O(M)$.
Traditionally the norm on a universal operator algebra of an
object is defined by taking the supremum over all representations
of the object as an operator algebra. Sometimes we are able to
define this norm without reference to the representations.   We
will call such a norm an intrinsic norm. We will construct an
intrinsic norm for the universal operator algebra of a directed
graph.

Let $M$ be a semigroup without identity.  It is well know
\cite{How:1976} that there is a monoid $M^+$ and a homomorphism
$\tau: M \rightarrow M^+$ which is one-to-one. It is a consequence
of the universal properties of the unitization \cite{Mey:2001} of
an operator algebra that there is a completely isometric inclusion
$\tilde{\tau} : O(M) \rightarrow O(M^+)$ induced by $ \tau$.
Further, by the definitions of $O(M^+)$ and $O(M)$ we know that
$\tilde{\tau}$ is a completely isometric isomorphism onto its
range, in particular $\| x \|_{O(M)} = \| \tau (x) \|_{O(M^+)}$.
The following lemma now follows.

\begin{lem} If $M$ is a semigroup without identity having no zero
divisors, then $O(M)^+ = O(M^+)$.\end{lem}

Now for $U \in M_n(\mathbb{C}M)$ define \[ \| U \|_n = \inf \{ \|
A_0\|\|A_1\| \cdots \|A_m \| \} \] where $A_0 \in M_{n,k}, A_m \in
M_{k,n} A_i \in M_k$ and $ U = A_0M_1A_2 \cdots M_mA_m$ where
$M_i$ is a diagonal $k\times k$ matrix with entries in $M^+$. It
is a consequence of \cite{Bl-Paul:1991} that if \[N = \{ x \in
\mathbb{C}M: \| x \|_1 = 0\] then $O(M)/N$ is an operator algebra
with the matricial norms given by $\| \cdot \|_n$.

We now turn to the context of universal operator algebras of
directed graphs.  We will denote by $W(Q)$ the set of all finite
words in the alphabet given by $E(Q) \cup V(Q)$.  We will denote
the range and source map by $r$ and $s$ respectively and we will
extend their definitions to all of $V(Q) \cup E(Q)$ be defining
$r(v) = s(v) = v$ for all $ v \in V(Q)$.

\begin{dfn} Let $Q$ be a directed graph and let $w \in W(Q)$.
We say that a word $ w_1w_2\cdots w_n \in W(Q)$ is \emph{reduced}
if $ r(w_i) \neq w_{i-1}$, $ s(w_i) \neq w_{i+1}$, $s(w_{i-1})
\neq w_i$, and $r(w_{i+1}) \neq w_i$. We denote by $w(Q)$, the set
of reduced words in $W(Q)$.\end{dfn}

\begin{prop} For a directed graph $Q$, $w(Q)$ is a semigroup.
Further $w(Q)$ has an identity if and only if $V(Q)$ is a
singleton.\end{prop}

\begin{proof} Certainly $W(Q)$ is a semigroup.  Further the
operation of reducing a word is terminating and locally confluent
and hence each word $w \in W(Q)$ has a unique reduced word $w(Q)$
associated to it.  It follows that $w(Q)$ is a semigroup, with
operation given by concatenation followed by reduction.

Now if $V(Q)$ has a single vertex $v$, then for an arbitrary edge
$e$ the reduced word for $ve $ is $e$, the reduced word for $e v$
is $e$ and the reduced word for $vv $ is $v$.  It follows that $v$
will serve as an identity element in $w(Q)$.

If $w(Q)$ has an identity element $\iota$ then $\iota^2 = \iota$
and hence $\iota$ corresponds to a vertex, since vertices give
rise to the only idempotents in $w(Q)$. On the other hand let $v$
be a vertex in $Q$. Then since $ \iota$ is an identity $\iota v =
v$ and hence $i= v$.
\end{proof}

We let $\|| \cdot \||_n $ denote the matricial norm on $O(w(Q))$.
It is a consequence of the proof of Proposition 2.1 in
\cite{Duncan:2004} that the subspace $ \{ x : \|| x \||_1 = 0 \}$
is the trivial subspace and hence $ \|| \cdot \||_n$ yields a norm
on $\mathbb{C}w(Q)$.

\begin{thm} Let $Q$ be a directed graph, then $O(w(Q))$ is
completely isometrically isomorphic to $OA(Q)$.\end{thm}

\begin{proof} Recalling the construction of $OA(Q)$ \cite{Duncan:2004} the
algebra $\mathbb{C}w(Q)$ is a dense subalgebra of $OA(Q)$. Further
the norm on $\mathbb{C}w(Q)$ is the universal norm induced by
representations of $\mathbb{C}w(Q)$.  The result now
follows.\end{proof}

It follows that matricial norms can be defined on $OA(Q)$ in an
intrinsic manner by defining matricial norms on $\mathbb{C}w(Q)$
as for the semigroup operator algebra.  This provides an intrinsic
characterization of the norm on $OA(Q)$ and perhaps makes the
construction of $OA(Q)$ less mysterious.

\begin{example} Let $T$ be the directed graph with two vertices
and a single edge connecting the vertices.  Labelling the vertices
as $v_0$ and $v_1$ and the edge as $t$, with $r(t) = v_1$ and
$s(t) = v_0$. we can see that $OA(T)$ is the norm closed algebra
generated by the span of elements of the set \[\{
(v_0)^{\delta_0}(v_1v_0)^{l_1}t^{n_1}(v_1v_0)^{m_1}(v_1v_0)^{l_2}t^{n_2}(v_1v_0)^{m_2}
\cdots (v_1v_0)^{l_k}t^{n_k}(v_1v_0)^{m_k}(v_1)^{\delta_1} \}\]
where $\delta_0, \delta_1 \in \{ 0, 1 \}$ and $ l_i, m_i, n_i \geq
0$.  This provides an alternate method from \cite{Duncan:2004}
where this algebra is described as the quotient of three free
products.\end{example}

\section{Intrinsic norms for $C^*_m(A)$}

We now look to build $C^*_m(A)$ in a manner intrinsic to $A$,
without reference to the completely contractive representations of
$A$.  This removes the need to define $C^*_m(A)$ by reference to
all completely contractive representations of $A$,as is done in
\cite{Bl:1999}.  In particular we will have a concrete
construction of the algebra $C^*_m(A)$ which should lead to a
better understanding of the maximal $C^*$-algebra of an operator
algebra.  We begin by letting $A$ be a unital operator algebra. We
recall the intrinsic characterization of the operator algebraic
free product of two operator algebras.

\begin{construction}[Blecher-Paulsen \cite{Bl-Paul:1991} Theorem 4.1] For $A$
and $B$ operator algebras with a common subalgebra $D$ and for $x$
in the algebraic free product of $A$ and $B$ amalgamated over $D$
we define
\[ \| x \|_{OA} = \inf \{ \| x_1 \| \| x_2 \| \cdots \| x_n \|
: x_1* x_2*\cdots * x_n = x \} \] where the $x_i$ are elements of
either $M_{k_i}(A)$ or $M_{j_i}(B)$ with $j_i, k_i \in
\mathbb{N}$, and $*$ the free product matrix multiplication.
Completing the algebraic free product with respect to this norm
yields an operator algebra $A*_{OA}B$ with the following universal
property.\end{construction}

Universal property:  If $\tau : A \rightarrow X$ and $\sigma : B
\rightarrow X$ are completely contractive such that $ \sigma|_D =
\tau|_D$ then there is $ \tau * \sigma : A*_{OA}B \rightarrow X$
completely contractive with $ \tau*\sigma|_A = \tau $ and $\tau *
\sigma |_B = \sigma$.

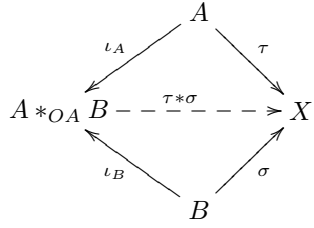
\[
\xymatrix{ & A \ar[dl]_{\iota_A} \ar[dr]^{\tau} & \\ A*_{OA}B
\ar@{-->}[rr]^{\tau * \sigma} & & X \\ & B \ar[ul]^{\iota_B}
\ar[ur]_{\sigma} & }
\]

Using this construction we will be able to build $C^*_m(A)$
intrinsically.  We begin with a definition.  Recall that the
diagonal of an operator algebra $A \cap A^*$ is independent of
representation and hence the diagonal can be found by taking any
faithful representation and finding the diagonal in that
particular representation.

\begin{dfn} Let $A$ be a unital operator algebra and let
$A^{\#}$ be the canonically associated adjoint algebra.  We write
$\mathcal{F(A)} = A*_{\Delta(A)}A^{\#}$ for the operator algebraic
free product amalgamated over the diagonal, $\Delta(A)$, of $A$.
\end{dfn}

\begin{rmk}  We will use ${\#}$ to denote the formal adjoint, and
$*$ to represent an adjoint in a $C^*$ algebra.  The reason is to
minimize confusion between elements of $\mathcal{F}$ and elements
of $C^*_m(A)$. \end{rmk}

\begin{rmk} Where it will not cause confusion we suppress
the $A$ in the notation that follows.\end{rmk}

We now construct a $C^*$ semi-norm on $\mathcal{F}$.

\begin{dfn} For $y \in M_n(\mathcal{F})$ we say that $y
\geq 0 $ if \[ y = \sum_{i=1}^m y_i^{\#}y_i \] for some set $ \{
y_i\}_{i=1}^m \in M_{{k_i},n}( \mathcal{F}),$ where $k_i$ is a
positive integer.
\end{dfn}

We record some elementary lemmas involving positive elements of
$\mathcal{F}$ and $M_2( \mathcal{F})$. In what follows we will
omit the formal product symbol where it is inferred.

\begin{lem} Let $x \in \mathcal{F}$ and $t \in \mathbb{R}^+$,
then \[ \begin{bmatrix} t & x \\ x^{\#} & t \end{bmatrix} \geq 0
\mbox{ if and only if } t^2 - x^{\#}x \geq 0. \] \end{lem}

\begin{proof} If $t^2 - x^{\#}x \geq 0 $ then clearly
\[ \begin{bmatrix} t^2 & 0 \\ 0 & t^2 - x^{\#}x \end{bmatrix}
\geq 0 \] and hence \[ \begin{bmatrix} t & x \\ x^{\#} & t
\end{bmatrix} = \frac{1}{\sqrt{t}} \begin{bmatrix} 1 & 0 \\
\frac{x^{\#}}{t} & 1\end{bmatrix}\begin{bmatrix} t & 0 \\
0 & t^2 - x^{\#}x\end{bmatrix} \begin{bmatrix} 1 & \frac{x}{t} \\
0 & 1\end{bmatrix}\frac{1}{\sqrt{t}} \geq 0. \]

Now if $\displaystyle{\begin{bmatrix} t & x \\ x^{\#} & t
\end{bmatrix} \geq 0} $ then \begin{eqnarray*} 0 & \leq &
\begin{bmatrix}1 & 0 \\ -\frac{x^{\#}}{t} & 1 \end{bmatrix} \begin{bmatrix}
t & x \\ x^{\#} & t\end{bmatrix} \begin{bmatrix}1 & -
\frac{x^{\#}}{t}
\\ 0 & 1 \end{bmatrix} \\ & = & \begin{bmatrix}1 & 0 \\ -\frac{x}{t}
& 1 \end{bmatrix} \begin{bmatrix}1 & 0 \\ \frac{x^{\#}}{t} & 0
\end{bmatrix} \begin{bmatrix} t & 0 \\ 0 & t -
\frac{x^{\#}x}{t}\end{bmatrix}
\begin{bmatrix} 1 & \frac{x}{t} \\ 0 & 1 \end{bmatrix} \begin{bmatrix}1 & - \frac{x}{t} \\
0 & 1 \end{bmatrix} \\ & = & \begin{bmatrix}t & 0 \\ 0 & t -
\frac{x^{\#}x}{t}\end{bmatrix}.
\end{eqnarray*}  It follows that $ t^2 - x^{\#}x \geq 0 $.
\end{proof}

The next two lemmas will allow us to show that the semi-norm we
define later is actually a $C^*$ semi-norm.

\begin{lem} For $x \in \mathcal{F}$ and $ t \in \mathbb{R}^+$ \[
\begin{bmatrix} t & x \\ x^{\#} & t \end{bmatrix} \geq 0 \mbox{ if and only if }
\begin{bmatrix} t^2 & x^{\#}x \\ x^{\#}x & t^2 \end{bmatrix} \geq 0. \] \end{lem}

\begin{proof}  If $ \displaystyle{\begin{bmatrix} t^2 & x^{\#}x \\ x^{\#}x & t^2 \end{bmatrix} \geq
0}$ then it follows that \begin{eqnarray*} 0 & \leq &
\begin{bmatrix} 1 & -1 \\ 0 & 0 \end{bmatrix} \begin{bmatrix} t^2 & x^{\#}x \\
x^{\#}x & t^2 \end{bmatrix} \begin{bmatrix}1 & 0
\\ -1 & 0 \end{bmatrix} \\ & = & \begin{bmatrix} 2t^2- 2x^{\#}x & 0 \\ 0 & 0 \end{bmatrix}
\end{eqnarray*} which implies that $t^2 - x^{\#}x \geq 0$ which by
the previous lemma yields one direction of the result.

Now if $\displaystyle{ \begin{bmatrix} t & x \\ x^{\#} & t
\end{bmatrix} \geq 0 }$ then by the previous lemma we have $ t^2 -
x^{\#}x \geq 0 $ and clearly $ t^2 + x^{\#}x \geq 0 $.  It follows
that $ \displaystyle{ \begin{bmatrix} t^2 + x^{\#}x & 0 \\ 0 & t^2
- x^{\#}x\end{bmatrix} \geq 0 } $ and hence
\begin{eqnarray*} 0 & \leq & \begin{bmatrix} 1 & -1 \\ 1 & 1\end{bmatrix}\begin{bmatrix}
t^2 + x^{\#}x & 0 \\ 0 & t^2 - x^{\#}x\end{bmatrix}\begin{bmatrix} 1 & 1 \\ -1 & 1\end{bmatrix} \\
& = & 2 \begin{bmatrix} t^2 & x^{\#}x \\ x^{\#}x & t^2
\end{bmatrix}.\end{eqnarray*}
\end{proof}

\begin{lem} If $x, y \in \mathcal{F}$, $s, t \in
\mathbb{R}^+$, \[ \begin{bmatrix}s & x \\ x^{\#} & s\end{bmatrix}
\geq 0, \mbox{ and }
\begin{bmatrix}t & y \\ y^{\#} & t\end{bmatrix} \geq 0 \] then \[ \begin{bmatrix} st & xy \\
y^{\#}x^{\#} & st \end{bmatrix} \geq 0. \] \end{lem}

\begin{proof} Notice by the first lemma that $s^2- x^{\#}x \geq 0 $
and hence $ s^2 y^{\#}y - y^{\#}x^{\#}xy \geq 0 $.  But notice
also that $ t^2 - y^{\#}y \geq 0 $ hence $ s^2t^2 - s^2 y^{\#}y
\geq 0 $ and it follows that $ s^2t^2 - y^{\#}x^{\#}xy \geq 0 $
and the first lemma gives the result. \end{proof}

The next lemma follows trivially from the definition.  The next
four are the last steps in providing a $C^*$-seminorm on
$\mathcal{F}$.

\begin{lem} Let $x, y \in \mathcal{F}$ and $ s, t \in
\mathbb{R}^+$ with \[ \begin{bmatrix} s & x \\ x^{\#} &
s\end{bmatrix} \geq 0 \mbox{ and } \begin{bmatrix} t & y \\ y^{\#}
& t \end{bmatrix} \geq 0 \] then
\[ \begin{bmatrix} s+t & x+y \\ x^{\#} + y^{\#} & s+t\end{bmatrix} \geq 0.\] \end{lem}

\begin{lem} Let $x \in \mathcal{F}$ and $ s \in \mathbb{R}^+ $
then \[ \begin{bmatrix} s & x \\ x^{\#} & s\end{bmatrix} \geq 0
\mbox{ if and only if } \begin{bmatrix} s & x^{\#} \\ x &
s\end{bmatrix} \geq 0. \]
\end{lem}

\begin{proof} Notice that \[ \begin{bmatrix} s & x \\ x^{\#} & s\end{bmatrix} =
\begin{bmatrix} 0 & 1 \\ 1 & 0 \end{bmatrix} \begin{bmatrix} s & x^{\#} \\ x & s\end{bmatrix}
\begin{bmatrix} 0 & 1 \\ 1 & 0\end{bmatrix} \] and the result is immediate.
\end{proof}

\begin{lem} Let $x \in \mathcal{F} $ and $ \lambda \in
\mathbb{C}, s \in \mathbb{R}^+$ then \[ \begin{bmatrix} | \lambda
| s & \lambda x \\ \overline{\lambda} x^{\#} & | \lambda | s
\end{bmatrix} \geq 0 \mbox{ if and only if } \begin{bmatrix} s & x
\\ x^{\#} & s \end{bmatrix} \geq 0.\]
\end{lem}

\begin{proof} Notice first that $ | \lambda |^2
s^2 - | \lambda |^2 x^{\#}x \geq 0$ if and only if $ s^2 - x^{\#}x
\geq 0 $. Now by the first lemma the result is established.
\end{proof}

\begin{lem} Let $x \in \mathcal{F}$ then \[ \inf \left\{ t:
\begin{bmatrix} t & x \\ x^{\#} & t \end{bmatrix} \geq 0 \right\} \leq \| x
\|_{OA}.\]
\end{lem}

\begin{proof} Notice that this is equivalent to showing that $ \|
x \|^2 - x^{\#}x \geq 0 $ in $\mathcal{F}$.  This will follow by
an induction.  We can begin by letting $x_1x_2$ be a matrix
factorization of $x$, and we define $t_i = \| x_i \| $ in the
appropriate matrix algebra.  Then notice that \begin{align*} t_1^2
t_2^2 - x^{\#} x & = t_2^2 - x_2^{\#} x_1^{\#} x_1 x_2 \\ &=
x_2^{\#}(t_1^2 - x_1^{\#} x_1) x_2 + t_1^2( t_2^2 - x_2^{\#} x_2
\end{align*} which is sum of positive elements of and hence \[ t_2
^2 t_1^2 \geq \inf \left\{ t: \begin{bmatrix} t & x \\ x^{\#} & t
\end{bmatrix} \right\}. \]  Repeating the process for larger
matrix factorization tells us that for any factorization $ x =
x_1x_2, \cdots x_n$ we have \[ \inf \left\{ t: \begin{bmatrix} t &
x \\ x^{\#} & t \end{bmatrix} \right\} \geq \| t_1\| \| t_2 \|
\cdots \| t_n \|. \] As the factorization is arbitrary the result
follows.\end{proof}

We are now in a position to define a $C^*$ seminorm on
$\mathcal{F}$.

\begin{dfn} Let $ x \in \mathcal{F}$ then define \[ \gamma (x)
= \inf \left\{ t: \begin{bmatrix} t & x \\ x^{\#} & t
\end{bmatrix} \geq 0 \right\}.
\]
\end{dfn}

\begin{prop}\label{intrinsic} The function $ \gamma$ is a
$C^*$ seminorm on $\mathcal{F}$ and $ \mathcal{F} / \ker \gamma $
is isomorphic to $ C^*_m(A)$.\end{prop}

\begin{proof}  The statement that $ \gamma$ is a $C^*$ seminorm
follows from the series of lemmas preceding the definition.  We
need only show that $\mathcal{F} / \ker(\gamma) \cong C^*_m(A)$.
Notice that $ q : A \rightarrow \mathcal{F} / \ker ( \gamma)$ is
completely contractive and hence the induced map $ q : \mathcal{F}
\rightarrow \mathcal{F} / \ker (\gamma)$ is a completely
contractive quotient homomorphism. Notice that the inclusion $
\iota : A \rightarrow \mathcal{F}$ is completely isometric, and
further $ q \circ \iota : A \rightarrow \mathcal{F} /\ker (\gamma)
$ sends $A$ to a generating subalgebra of $\mathcal{F} / \ker
(\gamma)$.  It follows by the universal property for $ C^*_m(A) $
that there exists an onto $*$ homomorphism $\tilde{q} : C^*_m(A)
\rightarrow \mathcal{F} / \ker ( \gamma)$.

We also know that there exists a completely contractive
homomorphism $\pi : \mathcal{F} \rightarrow C^*_m(A)$.  Now if $
\varepsilon - x^{\#}x \geq 0 $ for all $ \varepsilon > 0 $, then $
\pi( \varepsilon - x^{\#}x) \geq 0 $ for all $ \varepsilon > 0$.
In particular $ \varepsilon - \pi ( x)^* \pi(x) \geq 0 $ for all $
\varepsilon > 0$.  Now as $ C^*_m(A)$ is a $C^*$ algebra it
follows that $ \pi(x) = 0$.  It follows that $ \ker ( \gamma)
\subseteq \ker (\pi)$.  Hence there is a completely contractive
homomorphism $ \overline{\pi} : \mathcal{F} / \ker( \gamma)
\rightarrow C^*_m(A)$.  Notice that $ \overline{ \pi } \circ
\tilde{q} (x) = x$ for all $ x \in A$ it follows that $C^*_m(A)
\cong \mathcal{F} / \ker (\gamma) $.
\end{proof}

This would seem to imply the Blecher-Ruan-Sinclair Theorem (BRS
Theorem), see \cite[corollary 16.7]{Paul:2002}. Recall though that
the BRS-theorem was implicit in constructing the operator
algebraic free product \cite{Bl-Paul:1991} and hence this does not
provide an alternate approach to the BRS-theorem.

Notice that in \cite{Palm:2001} the function defined above is
defined for a general Banach $*$-algebra.  There it is shown that
the quotient is the maximal $C^*$-algebra representation of the
Banach $*$-algebra.  In particular we can use the general theory
of Banach-$*$ algebras to get at the same result \cite[proposition
11.1.4]{Palm:2001}. We need only show that $\mathcal{F}$ is indeed
a Banach-$*$ algebra.

\begin{prop} Let $A$ be an operator algebra then $
\mathcal{F} = A *_{\Delta(A)} A^{\#} $ is a Banach $*$-algebra.
\end{prop}

\begin{proof} We know that $\mathcal{F}(A)$ is an operator algebra and hence a
Banach algebra.  Now let $ x $ be in $ A *_{alg} A^{\#} $ and let
$ \varepsilon > 0 $.  By definition there exists  $A_1 \in
M_{k_1}(A) ,A_2 \in M_{k_2}(A) , \cdots , A_n \in M_{k_n}(A)$ and
$  B_1 \in M_{j_1}(A^{\#}), B_2 M_{j_2}(A^{\#}), \cdots B_n \in
M_{j_n}(A^{\#})$ such that  $ A_1* B_1* A_2* B_2* \cdots * A_n*
B_n = x$ and \[ \| x \|_{OA} \leq \| A_1 \|_{M_{k_1}(A)} \| B_1
\|_{M_{j_1}(A^{\#})} \cdots \| B_n \|_{M_{j_n}(A^{\#})} \leq \|x
\|_{OA} + \varepsilon. \] Now notice that $ x^{\#} = B_n^{\#}
\cdots A_1^{\#} $ and \begin{align*} \| A_1 \|_{M_{k_1}(A)} & \|
B_1 \|_{M_{j_1}(A^{\#})} \cdots \| B_n \|_{M_{j_n}(A^{\#})} \\ &=
\| B_n \|_{M_{j_n}(A)} \| A_n \|_{M_{k_n}(A^{\#})} \cdots \| A_n
\|_{M_{k_n}(A^{\#})} .\end{align*}  It follows that $ \| x^{\#}
\|_{OA} \leq \| x \|_{OA} $.  A similar argument tells us that $
\| x ^{\#} \|_{OA} \leq \| x \|_{OA} $ and hence $ \| x ^{\#}
\|_{OA} = \| x \|_{OA} $

Now if $ \{ x_n \}$ is Cauchy, then $ \{ x_n^{\#} \}$ is cauchy
and hence convergent.  Now if $ x_n \rightarrow x$ then $ \| x_n -
x \|_{OA} \rightarrow 0$.  By uniqueness of limits it follows that
$ \| x_n^{\#} - x^{\#} \|_{OA}  \rightarrow 0$. Hence $ \| \cdot
\|_{OA}$ is continuous with respect to $ {\#}$ and hence $
\mathcal{F}$ is a Banach $*$ algebra.\end{proof}

In \cite{Palm:2001} $ \gamma$ is called the Gelfand-Naimark
seminorm and the ideal $ \ker \gamma $ is called the reducing
ideal of $\mathcal{F}$.

In this section we have constructed $C^*_m(A)$ intrinsically for
an operator algebra $A$.  This construction, in particular,
applies to the algebra $GC^*_m(Q)$ where $Q$ is a directed graph.
Given a directed graph $Q$, we can use Theorem 3.3 in
\cite{Duncan:2004} to recognize $GC^*_m(Q)$ as a maximal $C^*$
envelope of $OA(Q)$. Proposition \ref{intrinsic} then gives us an
intrinsic seminorm on $OA(Q) * OA(Q)^*$ which yields the algebra
$GC^*_m(Q)$.

\section{Idempotents in $OA(Q)$ and $GC^*_m(Q)$}\label{ktheory}

We remind the reader of an example from \cite{Duncan:2004} and a
result concerning the $K$-groups of $OA(Q)$ and $GC^*_m(Q)$.

\begin{example}\label{dots} We will denote by $V_n$ the graph with
$n$ vertices and no edges.  $OA(V_n)$ is equal to the
unamalgamated free product of copies of $\mathbb{C}$.
\end{example}

\begin{prop} Let $Q$ be a directed graph, then there is norm
continuous homotopy from $OA(Q)$ onto $OA(V(Q))$ and also from
$GC^*_m(Q)$ onto $OA(V(Q))$.\end{prop}

\begin{cor} Let $Q$ be a finite directed graph.
Then \[ K_0(OA(Q)) = K_0(GC^*_m(Q)) = \mathbb{Z}^{|V(Q)|}.
\]\end{cor}

\begin{proof} The previous proposition tells us that \[ K_0(OA(Q)) = K_0(V(Q))
=K_0(GC^*_m(Q)).\]  By applying a result of a result of Cuntz
\cite{Cu:1982} to the algebra $OA(V(Q))$ we get $K_0(OA(V(Q)) =
K_0 (\mathbb{C})^{|V(Q)|} = \mathbb{Z}^{|V(Q)|}$. \end{proof}

It follows that the $K_0$-groups count the number of vertices.
Hence the number of vertices is a Banach algebra invariant of the
algebra.  We will see in the next two sections that more is true.
The maximal ideal space will allow us to not only count the
vertices but it will also be used to identify the projections $\{
P_v: v \in V(Q) \}$.

\section{The maximal ideal space of $OA(Q)$ and $GC^*_m(Q)$
for finite graphs}\label{maximalideal}

For the remainder of this chapter we will only be concerned with
finite graphs. For a Banach algebra $A$, we denote the maximal
ideal space by $M_A$.  By $\mathcal{P}(X)$ we mean the power set
of $X$ and we let $ \mathbb{P}(X) = \mathcal{P}(X) \setminus \{
\emptyset \}$. For $ k \in \mathbb{N} $ we let $
\overline{\mathbb{D}}^k$ be the cartesian product of $k$ copies of
$\mathbb{D}$, and we let $ \overline{\mathbb{D}}^0 = \{ 0 \}$. If
$ S \subseteq V(Q)$ and $S \neq \emptyset$ we let
\[ \mathcal{E}(S):= \{e \in E(Q): r(e), s(e) \in S \}. \]  Lastly,
for $S$ a nonempty subset of $V(Q)$, we define $ n(S) =
|\mathcal{E}(S)|$.

\begin{prop}\label{mispace} The set $M_{OA(Q)}$ is homeomorphic to
\[\bigsqcup_{S \in \mathbb{P}(V(Q))}
\overline{\mathbb{D}}^{n(S)}.\]\end{prop}

\begin{proof} We begin by letting $ \varphi$ be a multiplicative
linear functional and fixing an enumeration of $E(Q)$. Now $
\varphi$ is uniquely determined by $ \varphi(P_v)$ and $
\varphi(T_e)$ where $ v \in V(Q)$ and $e \in E(Q)$.  It is clear
that $ \varphi(P_v) \in \{ 0,1 \}$ where $ v \in V(Q)$. Further we
have that \[ \| \varphi(T_e) \| \leq \| \varphi \| \| T_e \| \leq
1\] for all edges $ e \in E(Q)$ and hence $ \varphi (T_e) \in
\overline{\mathbb{D}}$.

Fix $ \varphi \in M_{OA(Q)}$ and let \[ S_{\varphi}:= \{ v \in
V(Q): \varphi(P_v) = 1 \}. \]  Then $ \varphi$ is determined with
a fixed ordering on $ \mathcal{E}(S_{\varphi})$ by the
$n(S_{\varphi})$-tuple
\[ [ \varphi(T_{e_{S_1}}), \varphi(T_{e_{S_2}}), \cdots
\varphi(T_{e_{S_{n(S)}}})]. \]  Thus the map $ \varphi \mapsto [
\varphi(T_{e_{S_1}}), \varphi(T_{e_{S_2}}), \cdots
\varphi(T_{e_{S_{n(S)}}})]$ gives a map of the maximal ideal space
into \[\bigsqcup_{S \in \mathbb{P}(V(Q))}
\overline{\mathbb{D}}^{n(S)}.\] We claim that this correspondence
is onto.  Uniqueness follows by definition.

To prove that the correspondence is onto, let $S \in
\mathcal{P}(V(Q))$ be nonempty and take $ \lambda \in
\overline{\mathbb{D}}^{n(S)}$. We define $ \varphi_{\lambda} : Q
\rightarrow \mathbb{C}$ by \[ \varphi_{\lambda} (v) =
\begin{cases} 1 & v \in S \\ 0 & \mbox{else}\end{cases}. \]  Then
define $ \varphi_{\lambda}(e_i) = \lambda_i$ for $ e_i \in \{ e :
r(e) \in S \mbox{ and } s(e) \in S \}$, and $ \varphi_{\lambda}(e)
= 0$ otherwise.  It is easy to see that $ \varphi_{\lambda}$ is a
contractive representation of $Q$. Now by the universal property
of $ OA(Q)$ there exists a unique completely contractive
homomorphism, which we also call $ \varphi_{\lambda}$, with $
\varphi_{\lambda}: OA(Q) \rightarrow \mathbb{C}$.  It follows that
the correspondence is onto.

We now turn to continuity.  If $\varphi_{\lambda} \rightarrow
\varphi$ then $ \varphi_{\lambda}(T_e) \rightarrow \varphi(T_e)$
and $ \varphi_{\lambda}(P_v)\rightarrow \varphi(P_v)$ for each
edge $e$ and vertex $v$.  It follows that the correspondence will
preserve the set $S$ and the $n(S)$ tuples will converge
pointwise. Thus the correspondence induces a continuous map
between $ M_{OA(Q)}$ and $\bigsqcup_{S \in \mathcal{P}(V(Q))}
\overline{\mathbb{D}}^{n(S)}$.  Now since we have a one to one and
onto continuous map from a space which is Hausdorff and compact we
have that the inverse map is also continuous and the homeomorphism
is established.\end{proof}

In fact we have established that $M_{OA(Q)}$ is a compact
Hausdorff space with a connected component for each nonempty $ S
\subseteq V(Q)$.

\begin{example} Let $Q$ be the graph 
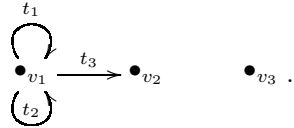
\[ \xymatrix{ {\bullet}_{v_1} \ar@(ul,ur)[]^{t_1}
\ar@(dl,dr)[]^{t_2} \ar[r]^{t_3} & {\bullet}_{v_2} &
{\bullet}_{v_3}}. \] Since $Q$ has 3 vertices there are seven
connected components in the maximal ideal space.  The component
corresponding to $v_1$ has two copies of $\overline{\mathbb{D}}$
since there are two edges with range and source equal to $v_1$.
The component corresponding to $v_2$ and the component
corresponding to $v_3$ are both singleton sets since neither
vertex has an edge which enters and leaves the vertex. The
component corresponding to the pair $\{ v_1, v_2 \}$ has three
copies of $\overline{\mathbb{D}}$ one for each of the edges,
$t_1,t_2,$ and $t_3$.  The pair $\{ v_2,v_3\}$ also yields a
singleton set.  The pair $\{ v_1, v_3 \}$ has two copies of
$\overline{\mathbb{D}}$.  The final component corresponding to $\{
v_1,v_2,v_3 \}$ has 3 copies of $\overline{\mathbb{D}}$ since
there are three total edges.  The maximal ideal space is then
homeomorphic to \[ \overline{\mathbb{D}}^2 \sqcup
\overline{\mathbb{D}}^0 \sqcup \overline{\mathbb{D}}^0 \sqcup
\overline{\mathbb{D}}^3 \sqcup \overline{\mathbb{D}}^0 \sqcup
\overline{\mathbb{D}}^2 \sqcup
\overline{\mathbb{D}}^3.\]\end{example}

\begin{dfn} For a finite directed graph $Q$ we let $N_Q$ be the
number of connected components of $M_{OA(Q)}$.\end{dfn}

We can actually define several invariants of the algebra by using
combinatorial arguments and the structure of the maximal ideal
space in a fairly simple manner.

\begin{prop} For a finite directed graph $Q$, \[|V(Q)| = \log_2 ( N_Q
+ 1 )\] and
\[ |E(Q)| = \max\{ n(S) \}\] where $n(S) = |\mathcal{E}(S)|$.\end{prop}

\begin{proof} Each connected component of $M_{OA(Q)}$ is associated
uniquely to a nonempty subset of $V(Q)$.  It follows that $N_Q + 1
= | \mathcal{P}(V(Q))| = 2^{|V(Q)|}$. and the first formula is
established.  Secondly as $V(Q) \in \mathbb{P}(V(Q))$ there is a
connected component of $M_{OA(Q)}$ associated to the set $V(Q)$.
But $n(V(Q))$ is the number of edges emanating from and ending in
$V(Q)$, which is the total number of edges. Since $n(S) $ is less
than or equal to the total number of edges for all $S \subseteq
V(Q)$, we have the second formula and the corollary is
established.\end{proof}

\begin{cor}\label{counting2} Suppose $OA(Q_1)$ and $OA(Q_2)$ are algebraically isomorphic.
Then $ |V(Q_1)| = |V(Q_2)|$, $|E(Q_1)| = |E(Q_2)|$.\end{cor}

\begin{proof}If the algebras $OA(Q_1)$ and $ OA(Q_2) $ are isomorphic, then the
spaces $M_{OA(Q_1)}$ and $ M_{OA(Q_2)}$ are homeomorphic. Hence by
the formulas established in the previous Proposition the corollary
follows.\end{proof}

Actually, more is true.  We say that an edge $e \in E(Q)$ is a
loop edge if $s(e) = r(e)$.  We can use calculations to find the
number of loop edges and non loop edges in the graph from
combinatorial facts about $M_{OA(Q)}$.

\begin{prop}\label{counting} Let $Q$ be a finite directed graph.  If $\alpha$ is the number of loop
edges in $Q$ and $ \beta$ is the number of non loop edges in $Q$
then $\alpha$ and $ \beta$ can be calculated uniquely from
$M_{OA(Q)}$.
\end{prop}

\begin{proof} If  $n$ is the number of vertices in $Q$, then for an edge $e$ there
will be a copy of $\overline{\mathbb{D}}$ for every subset of $S
\subseteq V(Q)$ with $ r(e), s(e) \in S$. Thus if $r(e) = s(e)$,
since there are $2^{n-1}$ nonempty subsets of $V(Q)$ containing
$r(e)$, there are $2^{n-1}$ copies of $ \overline{\mathbb{D}} \in
M_{OA(Q)}$ for each loop edge.  If $r(e) \neq s(e)$ there are
$2^{n-2}$ subsets of $\mathbb{P}(V(Q))$ which contain $s(e)$ and
$r(e)$.  Thus, there are $ 2^{n-2}$ copies of $
\overline{\mathbb{D}} \in M_{OA(Q)}$ for each edge which is not a
loop.  If $ \alpha $ is the number of loop edges, and $ \beta $ is
the number of non loop edges, then $ \alpha (2^{n-1} ) + \beta (
2^{n-2}) = \sum_{S \in \mathbb{P}(V(Q))} n(S)$.  Now assume that
there are $ \alpha '$ and $\beta '$, a different combination of
loop edges and non loop edges respectively, such that $\alpha
(2^{n-1} ) + \beta ( 2^{n-2}) = \alpha' (2^{n-1} ) + \beta' (
2^{n-2})$.  Since the number of edges is fixed at $n$ we know that
$ (n-\beta)(2^{n-1} ) + \beta ( 2^{n-2}) = (n - \beta') (2^{n-1} )
+ \beta' ( 2^{n-2}) $. Simplifying we get that $ \beta (2^{n-2}-
2^{n-1}) = \beta ' (2^{n-2}- 2^{n-1})$ and hence $ \beta = \beta
'$.  It follows that in a finite graph the number of loop edges
and the number of non loop edges is an isomorphism invariant which
can be calculated directly from information about the set
$M_{OA(Q)}$.\end{proof}

We now look at the algebra $GC^*_M(Q)$.  Recall that this is the
universal $C^*$ algebra of the directed graph $Q$ which is
constructed by looking at $*$ representations of the graph $Q
*_{V(Q)}Q^{\leftarrow}$.  The universal properties will allow us
to identify the maximal ideal space of $GC^*_M(Q)$.

\begin{prop}\label{maximal}  The set $M_{GC^*_M(Q)}$ is homeomorphic to $
M_{OA(Q)}$.  In fact for an operator algebra $A$, $M_A$ is
homeomorphic to $M_{C^*_m(A)}$.\end{prop}

\begin{proof}  Since $GC^*_M(Q) = C^*_m(OA(Q))$ we will prove the
more general result.  If $ \varphi: C^*_m(A) \rightarrow
\mathbb{C}$ is a multiplicative linear functional then $
\varphi|_{A} \rightarrow \mathbb{C}$ is also a multiplicative
linear functional.  Further every multiplicative linear functional
$ \pi: A \rightarrow \mathbb{C}$ is completely contractive and
hence there exists a unique multiplicative linear functional $
\widetilde{\pi} : C^*_m(A) \rightarrow \mathbb{C}$ such that $
\widetilde{\pi}|_A = \pi|_A$. It follows that there is a one to
one correspondence between maximal ideals of $A$ and $C^*_m(A)$.
That the maps are continuous is trivial.\end{proof}

Proposition \ref{counting} applies also to $GC^*_m(Q)$ and hence
the number of vertices, loop edges and non loop edges are
isomorphism invariants for $GC^*_m(Q)$.

\section{Uniqueness of $GC^*_m(Q)$ for finite
graphs}\label{uniquecmg}

In this section we are interested in uniqueness of $GC^*_m(Q)$. We
begin with definitions. Let $A$ be a Banach algebra, and let $
\varphi $ be a multiplicative linear functional.  We let
\[ P({\varphi}) := \{ x \in A: x^2 = x \mbox{ and } \varphi(x) = 1 \}. \]

\begin{dfn} For $X$ a connected component of $M_A$ we say that $X$
has \emph{degree 1} if for every $ \varphi \in X$, $|P({\varphi})|
= 1$. We say that $X$ has \emph{degree k } for $ k > 1$ if there
are exactly $k$ degree 1 components $X_j$ such that $P({\varphi})
\cap \P({\tau}) \neq \emptyset$ for all $ \varphi \in X, \tau \in
X_j$. \end{dfn}

In the context of universal graph operator algebras the preceding
definition will be useful in establishing uniqueness.  We use it
now to identify the number of vertices associated to a particular
connected component of $M_{GC^*_m(Q)}$.  Recall that to each set $
S \in \mathbb{P}(X)$ there is a connected component in $
M_{GC^*_m(Q)}$.

It is a consequence of Proposition \ref{mispace} and Proposition
\ref{maximal} that for $ \varphi , \tau \in X$, a connected
component of $M_{GC^*_m(Q)}$, $P(\varphi) = P(\tau)$.

\begin{prop}\label{degreeone}  Let $X$ be a connected component in $M_{GC^*_m(Q)}$,
then $X$ has degree $k \geq 1$ if and only if there is a set of
disjoint vertices $V:= \{ v_1, v_2, \cdots, v_k \} \subseteq V(Q)$
such that $ X$ is the component associated to the set $ V \in
\mathbb{P}(X)$.\end{prop}

\begin{proof} Let $S \in \mathbb{P}(V(Q))$. We will show that if
$ |S| = k$ then the associated component, $X$, has degree $k$. let
$\{ v_1, v_2, \cdots, v_k \} = S$ and denote by $P_{v_i}$ the
projection associated to $v_i$.  Now define a contractive
representation $\pi: Q \rightarrow \mathbb{C}$ by sending $v_i$ to
$1$ for each $i$ and everything else to 0.  The induced completely
contractive map will be a multiplicative linear functional
associated to the component $X$. Now notice that $ P(\varphi) = \{
P_{v_1}, P_{v_2}, \cdots, P_{v_k} \}$ and hence $|P({\varphi})| =
k$.  The result follows.\end{proof}

\begin{dfn} For $X$ a connected component of degree k in $M_{GC^*_m(Q)}$ we
let \[ P_X := \bigcup_{Y \mbox{ degree } 1} \left(\bigcap
_{\varphi \in X, \tau \in Y} (P(\varphi) \cap P(\tau)) \right).
\] If $X$ has degree two and $Y$ and $Z$ are components of degree
one, then we say that $X$ is the component \emph{associated to $Y$
and $Z$} if $P_{X} = P_Y \cup P_Z$\end{dfn}

We are now in a position to describe the main result of this
section. Starting with the graph $Q$ we build an associated
undirected graph $\widehat{Q}$.  Recall that an undirected graph
is a 3-tuple $(V, E, n)$, where $V$ is a set of vertices, $E$ is
the set of all pairs of vertices, and $n: E \rightarrow
\mathbb{N}$ is a continuous map.  The map $n(\{ v, w \})$ will
specify how many edges connect the pair of vartices $v$ and $w$.
For a directed graph $Q$ we let $V(\widehat{Q}) = V(Q)$ and $n(\{
v, w\})$ be the number of edges $e$  with $\{ r(e), s(e) \} = \{
v, w \}$.  The graph $ \widehat{Q}$ can be thought of as the graph
obtained from $Q$ by removing the directions on each edge. We will
show that $\widehat{Q}$ is an isomorphism invariant for
$GC^*_m(Q)$.

We say that an edge $e$ in a graph is a loop edge if $r(e) =
s(e)$.

\begin{thm} Let $Q_1$ and $Q_2$ be finite directed graphs.  The
algebras $GC^*_m(Q_1)$ and $GC^*_m(Q_2)$ are isomorphic as Banach
algebras if and only if the graphs $\widehat{Q_1} $ and
$\widehat{Q_2}$ are isomorphic.\end{thm}

\begin{proof}  We begin by assuming that two algebras $GC^*_m(Q_1)$
and $GC^*_m(Q_2)$ are isomorphic.  It is well known that if two
Banach algebras are isomorphic as Banach algebras then there is an
induced homeomorphism between their maximal ideal spaces. It
follows from Corollary \ref{counting2} that the number of vertices
in $Q_1$ is equal to the number of vertices in $Q_2$. The
homeomorphism will clearly preserve degree. Let $X$ be a degree
two component in $M_{GC^*_m(Q_1)}$ or $M_{GC^*_m(Q_2)}$.  We can,
by Proposition \ref{degreeone}, identify the degree one
components, $Y$ and $Z$, which correspond to $X$. For clarity of
presentation we will write a degree two component with
corresponding degree one components $Y$ and $Z$ as $X_{Y,Z}$.

For an arbitrary connected component $W$ of $M_{GC^*_m(Q)}$ or
$M_{GC^*_m(Q')}$ we let $n(W)$ be the number of copies of
$\overline{\mathbb{D}}$ in $W$.  If $Y$ is a degree one component
then $n(Y)$ is the number of loop edges.  For a degree two
component $n(X_{Y,Z}) - (n(Y) + n(Z))$ is the number of edges $e$
in the graph with $r(e) \neq s(e)$ and $\varphi(r(e)) =1 =
\varphi(s(e))$ for all $ \varphi \in X_{Y,Z}$.  Now $n(W)$ is also
invariant under isomorphism and hence $\widehat{Q_1} $ is
isomorphic to $\widehat{Q_2}$.

For the converse, assume that $\widehat{Q_1}$ and $\widehat{Q_2}$
are isomorphic.  Then there is a 1-1 correspondence between the
sets $V(Q_1)$ and $V(Q_2)$.  There is also a 1-1 correspondence
between the sets $E(Q_1)$ and $E(Q_2)$.  We build a new directed
graph $\overline{Q}$ from $ \widehat{Q}$ by setting
$V(\overline{Q}) = V(\widehat{Q})$, $E(\overline{Q}) =
E(\widehat{Q})$, $r((v,w)) = v,$ and $ s((v,w))=w$. The assignment
of range and source will not change the graph
$\overline{Q}*_{V(\overline{Q})}\overline{Q}^{\leftarrow}$. It is
easy to see that
$\overline{Q_1}*_{V(\overline{Q_1})}\overline{Q_1}^{\leftarrow}$
is isomorphic to
$\overline{Q_2}*_{V(\overline{Q_2})}\overline{Q_2}^{\leftarrow}$.
The result now follows from the construction of $GC^*_m(Q)$,
\cite{Duncan:2004}.\end{proof}

There is no reason to expect a stronger uniqueness result.  For
example, the two graphs
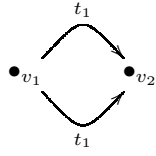
\[ \xymatrix{ {\bullet}_{v_1} \ar@(ur,ul)[r]^{t_1} \ar@(dr,dl)[r]_{t_1}
& {\bullet}_{v_2} } \] and 
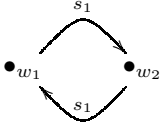
\[ \xymatrix{ {\bullet}_{w_1} \ar@(ur,
ul)[r]^{s_1} & {\bullet}_{w_2}  \ar@(dl,dr)[l]_{s_1}} \] have
isomorphic universal $C^*$ algebras, even though the graphs are
not isomorphic.  On the other hand any uniqueness may be
considered surprising since $Q$ is not an invariant for $C^*(Q)$.
For an example we point the reader to \cite{Hong-Szym:2003} where
it is shown that the distinct graphs 
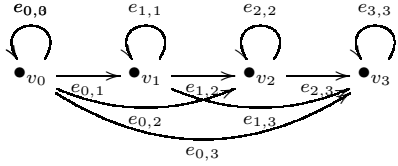
\[ \xymatrix{ {\bullet}_{v_0}
\ar@(ur,ul)[]_{e_{0,0}} \ar[r]_{e_{0,1}} \ar@/_1pc/[rr]_{e_{0,2}}
\ar@/_2pc/[rrr]_{e_{0,3}} \ar@(ur,ul)[]_{e_{0,3}} &
{\bullet}_{v_1} \ar@(ur,ul)[]_{e_{1,1}} \ar@/_1pc/[rr]_{e_{1,3}}
\ar[r]_{e_{1,2}} & {\bullet}_{v_2} \ar@(ur,ul)[]_{e_{2,2}}
\ar[r]_{e_{2,3}} & {\bullet}_{v_3} \ar@(ur,ul)[]_{e_{3,3}} } \]
and 
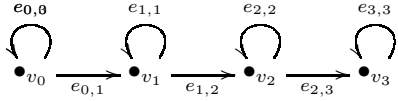
\[ \xymatrix{ {\bullet}_{v_0} \ar@(ur,ul)[]_{e_{0,0}}
\ar[r]_{e_{0,1}} \ar@(ur,ul)[]_{e_{0,3}} & {\bullet}_{v_1}
\ar@(ur,ul)[]_{e_{1,1}} \ar[r]_{e_{1,2}} & {\bullet}_{v_2}
\ar@(ur,ul)[]_{e_{2,2}} \ar[r]_{e_{2,3}} & {\bullet}_{v_3}
\ar@(ur,ul)[]_{e_{3,3}} } \] yield isomorphic $C^*$-algebras.

\section{Uniqueness of $OA(Q)$ for finite graphs}

We use the definitions from the previous section in establishing
the uniqueness of the algebra $OA(Q)$.  Once again we have most of
the information about our directed graph embedded in the maximal
ideal space. The only complication that remains is identifying the
directions on the edges with distinct source and range. We will
use ideas similar to those in \cite{Kat-Kribs:2003} to build the
original graph $Q$ from information about a class of
representations of $OA(Q)$.

As in Section \ref{uniquecmg} we will need to identify the degree
one components of $M_{OA(Q)}$.  The same arguments work in this
context so we do not repeat them here.  We need a few preliminary
results before addressing uniqueness.

Let $T_2$ be the algebra of $2 \times 2$ upper triangular
matrices.  If $A$ is an operator algebra we say that a
representation $\pi: A \rightarrow T_2$ is a \emph{two dimensional
nest representation} if $\pi$ is onto. Let $X$ and $Y$ be degree
one connected components in $M_A$. We say that a two dimensional
nest representation $\pi$ of $OA(Q)$ has the \emph{projection
property for $X$ and $Y$} if
\[ \pi (P_X) =
\begin{bmatrix} 1 & 0 \\ 0 & 0 \end{bmatrix}, \]\[ \pi
(P_Y) =
\begin{bmatrix} 0 & 0 \\ 0 & 1 \end{bmatrix}, \] and if $x$ is an
idempotent not contained in $P_X \cup P_Y$ then
\[ \pi (x) = 0.\]

\begin{dfn} For $X$ and $Y$, degree one connected components in
$M_{OA(Q)}$, let  \[K_{X,Y} := \cap \{ \ker(\pi): \pi \mbox{ has
the projection property for $X$ and $Y$} \} . \] \end{dfn}

If $\pi: OA(Q) \rightarrow T_2$ has the projection property for
$X$ and $Y$, then there is an induced map $\pi_q: OA(Q)/K_{X,Y}
\rightarrow T_2$ which is a two dimensional nest representation
with the projection property for $X$ and $Y$.

\begin{dfn} If $X$ and $Y$ are degree one connected components in
$M_{OA(Q)}$ let $R_{X,Y}$ be the set of all cosets $OA(Q) +K_{X,Y}
\in OA(Q)/K_{X,Y}$ such that \[(\pi_q(OA(Q)+K_{X,Y}))^2= 0\] for
all $\pi$ with the projection property for $X$ and $Y$.\end{dfn}

\begin{lem} $R_{X,Y}$ is a closed two sided ideal in
$OA(Q)/K_{X,Y}$. \end{lem}

\begin{proof} The fact that $(\pi_q(OA(Q)+K_{X,Y}))^2 = 0$ implies
that $\pi_q(OA(Q)+K_{X,Y})$ is strictly upper triangular.  Now if
$B+K_{X,Y}$ is another coset in $OA(Q) / K_{X,Y}$ then
\begin{align*}
\pi_q((A+K_{X,Y} )(B+K_{X,Y})) & = \pi_q(A+K_{X,Y}) \pi_q(B+K_{X,Y}) \\
& = \begin{bmatrix}0 & a \\ 0 & 0 \end{bmatrix} \begin{bmatrix}
a_1
& a_2 \\ 0 & a_3 \end{bmatrix} \\ &= \begin{bmatrix} 0 & a a_3 \\
0 & 0 \end{bmatrix}. \end{align*} Hence $(\pi_q(( A+K_{X,Y})
(B+K_{X,Y}) ))^2 = 0$.  Similar arguments for multiplication on
the right by an ideal element shows that $R_{X,Y}$ is a two sided
ideal. Closure is automatic since $\pi_q$ is
continuous.\end{proof}

We now describe $R_{X,Y}$ for a finite directed graph $Q$.

\begin{prop} Let $Q$ be a finite graph, and $v,w \in V(Q)$. We denote
the connected components of $M_{OA(Q)}$ associated to $\{ v \}$
and $\{ w \}$ by $V$ and $W$, respectively. There are $n$ edges
with range $v$ and source $w$ if and only if $R_{V,W}$ has a
minimal generating set of cardinality $n$.\end{prop}

\begin{proof}  Let $V$ and $W$ be the sets described.  Then it
is clear that $P_v A P_w + K_{V,W} = A + K_{V,W}$ for all $A \in
R_{V,W}$.  Further notice that $T_e + K_{V,W} = P_vT_eP_w +
K_{V,W} \in R_{V,W}$ for all edges $e$ with $s(e) = w$, $r(e) =
v$.  A quick calculation tells us that if $\pi$ has the projection
property for $V$ and $W$ then for each edge $e$ with $r(e) = v,
s(e) = w$ there is an $a_e$ such that
\[  \pi(T_e) = \begin{pmatrix}0 & a_e \\
0 & 0 \end{pmatrix}. \]  If $r(e) = s(e) = v$ then \[ \pi(T_e) =
\begin{bmatrix} \lambda & 0 \\ 0 & 0 \end{bmatrix} \] for some
$\lambda \in \mathbb{C}$.  Similarly if $r(e) = s(e) = w$ then \[
\pi(T_e) =
\begin{bmatrix} 0 & 0 \\ 0 & \lambda \end{bmatrix} \] for some
$\lambda \in \mathbb{C}$. Lastly if $r(e) \neq v$ or $s(e) \neq w$
then $ \pi(T_e) = 0$.  Now, letting \[ N = \{ e: r(e) = v, s(e) =
w \}\] we have that $\{ T_e + K_{V,W} : e \in N \}$ is a linearly
independent generating set for $R_{V,W}$.  In particular, a
typical element of $R_{V,W}$ is contained in the closure of the
linear span of
\[R:=\{ T_f^n T_e T_g^m + K_{V,W}: s(f) = r(f) = r(e), e \in N,
r(g) = s(g) = s(e), n,m \geq 0 \}.\]

Now let $X = \{ x_1, x_2, \cdots,x_m \}$ be a generating set for
$R_{V,W}$.  We will show that $|X| \geq n$.  Let $e$ be an an edge
in $N$.  Now define a representation $\pi_e: Q \rightarrow T_2$ by
first letting
\begin{align*} \pi_e({s(e)}) &=
\begin{bmatrix} 0 & 0 \\ 0 & 1 \end{bmatrix} \\ \pi_e({r(e)}) &=
\begin{bmatrix} 1 & 0 \\ 0 & 0 \end{bmatrix} \\ \pi_e(e) &=
\begin{bmatrix} 0 & 1 \\ 0 & 0 \end{bmatrix}\end{align*} and
sending all other edges and vertices to zero.  There will be a
completely contractive extension $\tilde{\pi_e}: OA(Q) \rightarrow
T_2$ and $\tilde{\pi_e}$ will have the projection property for $V$
and $W$.  It follows then that $ (\tilde{\pi_e})_q: OA(Q) /
K_{V,W} \rightarrow T_2$ is well defined.  Now there exists $x_i
\in X$ such that $ \| (\tilde{\pi_e})_q (x_i) \| > 0$ hence $ x_i
= \alpha_e T_e + k_e$ where $\alpha_e \neq 0 \in \mathbb{C}$ and
$k_e$ is in the kernel of $\tilde{\pi_e}$.

It follows that for each $i$, $x_i = \{ (\sum_{e \in N}
\alpha_eT_e) + k_i \}$ where $k_i \in \cap_{e \in N} \ker(\pi_e)$.
The set $\{ T_e \}_{e \in N} + \cap_{e \in N} \ker(\pi_e)$ is a
linearly independent subset of $OA(Q) / \cap_{e \in N}
\ker(\pi_e)$ and it follows that $|X| > n$.\end{proof}

We now prove a classification theorem for universal operator
algebras of directed graphs.

\begin{thm} Let $Q_1$ and $Q_2$ be finite directed graphs. The
algebras $OA(Q_1)$ and $OA(Q_2)$ are isomorphic as Banach algebras
if and only if $Q_1$ and $Q_2$ are isomorphic as directed
graphs.\end{thm}

\begin{proof} Certainly if $Q_1$ and $Q_2$ are isomorphic then
$OA(Q_1)$ and $OA(Q_2)$ are isomorphic by uniqueness of the
extension of a directed graph morphism.  Let the map $\pi:
OA(Q_1)\rightarrow OA(Q_2)$ be a bounded isomorphism then $\pi$
induces a homeomorphism, which we also call $\pi$. Further the
algebra $OA(Q_1)/K_{V,W}$ is isomorphic to $OA(Q_2)/ K_{\pi(V),
\pi(W)}$, and the result follows.
\end{proof}

We note here some differences between the proof here and the proof
of uniqueness for quiver algebras given in \cite{Kat-Kribs:2003}.
For the quiver algebras, since the projections are orthogonal, all
of the connected components of the maximal ideal space are degree
one.  This simplifies the quiver algebra result.  Also, Katsoulis
and Kribs, use a Fourier expansion for elements of the quiver
algebra and hence they do not need to restrict the class of two
dimensional nest representations that they use. Although the
proofs are significantly different the ideas are similar. Using
the maximal ideal space and the two dimensional nest
representations we construct the underlying graph from the algebra
using only Banach algebra properties.

Extending the uniqueness results to infinite graphs may be more
complicated. A better understanding of the maximal ideal space is
vital.  On the other hand it is a consequence of the description
of the $K$-theory that if $Q$ is an infinite graph, then $OA(Q)
\not\cong OA(Q')$ for any finite graph $Q'$ and similarly
$GC^*_m(Q) \neq\cong OA(Q')$ for any finite graph $Q'$.

\begin{acknowledgement}  Part of this work is contained in the authors
doctoral dissertation.  He gratefully acknowledges the support of
his department and advisor, David Pitts.\end{acknowledgement}

\bibliographystyle{plain}
\bibliography{common}

\end{document}